\documentclass{article}
\usepackage{amsfonts}
\usepackage[ruled,lined]{algorithm2e}
\usepackage{amsmath}		
\usepackage{amssymb}
\usepackage{cite}
\usepackage{float}
\usepackage{graphicx}
\usepackage{booktabs, multicol, multirow}
\usepackage{bigstrut}
\usepackage{array}
\usepackage{tabularx}
\usepackage{pdflscape}
\usepackage{adjustbox}
\usepackage{longtable}
\usepackage{makeidx}
\usepackage{rotating}
\usepackage{wasysym}
\usepackage{graphicx,xcolor}
\usepackage{graphics}
\usepackage{float}
\usepackage{latexsym}
\usepackage{amsfonts}
\usepackage{mathptmx}
\usepackage{makeidx}
\usepackage[latin2]{inputenc}
\usepackage{dsfont}
\usepackage[sc]{mathpazo} 
\usepackage[T1]{fontenc} 
\usepackage{microtype} 

\usepackage[hmarginratio=1:1,top=23mm,columnsep=20pt]{geometry} 
\usepackage{chngpage} 
\usepackage{hyperref} 

\usepackage{booktabs} 
\usepackage{float} 

\usepackage{lettrine} 
\usepackage{paralist} 

\usepackage{abstract} 
\usepackage{titlesec} 
\titleformat{\section}[block]{\large\scshape\centering}{\thesection.}{1em}{} 
\usepackage{fancyhdr} 
\newtheorem{theorem}{Theorem}[section]
\newtheorem{lemma}[theorem]{Lemma}

\newtheorem{remark}[theorem]{Remark}

\newtheorem{assumption}[theorem]{Assumption}
\def\QED{~\rule[-1pt]{6pt}{6pt}\par\medskip}

\newenvironment{proof}{{\bf Proof\ }}{ \hfill \QED}

\begin{document}
\title{\vspace{-15mm}\fontsize{16pt}{10pt}\selectfont{A multivariate spectral hybridization of HS and PRP
method for nonlinear systems of equations
}} 
\author{Hassan Mohammad$^{\rm 1}$\\ [0.4em]  
  {\footnotesize {\it  $^{\rm 1}$ Numerical Optimization Research Group, Department of Mathematical Sciences}},\\ {\footnotesize {Faculty of Physical Sciences, Bayero University, Kano, Nigeria.}}
  \\ [-0.4em]
}
 
 \date{}    
\maketitle 



\begin{abstract}
\noindent We present a multivariate spectral hybridization of   Hestenes-Stiefel (HS) and Polak--Ribi{\`e}re--Polyak (PRP) method for solving large-scale nonlinear systems of equations. The search direction of the method is obtained by incorporating a multivariate spectral approach with the positive hybridization of Hestenes-Stiefel and Polak--Ribi{\`e}re--Polyak parameters (HS \& PRP hybrid+). By employing a derivative-free nonmonotone line search technique, the global convergence of the sequence generated by the method is proven. Numerical experiments are given to demonstrate the good performance of the method compared with similar methods in the literature designed for solving large-scale nonlinear systems of equations.
\end{abstract}
\textbf{Keywords}: Nonlinear equations, derivative-free methods, multivariate spectral method, hybridization approach, global convergence, numerical results \\\\
\textbf{Mathematics Subject Classification}: 65H10, 65H20, 90C06, 90C52, 90C56.
\section{Introduction}
\label{intro}
We consider the problem of solving nonlinear equations of the form
\begin{equation}\label{eq1}
	F(x)=0, 
\end{equation}
where $F : \mathbb{R}^n \to \mathbb{R}^n$ is a continuously differentiable function. The solution of a system of nonlinear equations is one of the important tasks in numerical analysis and optimization. Many applied problems from sciences, engineering and economics can be reduced to solving nonlinear equations. Newton's method is one of the well-known methods for solving \eqref{eq1}, however, because of the requirement of computing the Jacobian matrix which requires first-order derivatives, the application of Newton's method is limited especially for nonsmooth problems. For solving large-scale systems of nonlinear equations, matrix-free approaches such as spectral methods \cite{zhang2006,lacruz2006,mohd17,huang2017}, conjugate residual methods \cite{li2011,li2014derivative,liu2018d,cheng2009prp} are preferable to Newton's like methods because they are simple to implement and require less amount of storage.

Recalling the unconstrained optimization problem
\begin{equation}\label{eq2}
	\min _{x\in \mathbb{R}^n}f(x),
\end{equation}
where $f:\mathbb{R}^n\to \mathbb{R}$ is a continuously differentiable function with gradient $g(x)=\nabla f(x)$.   

Conjugate gradient (CG) methods are class of matrix-free methods that are efficient with good global convergence properties for solving \eqref{eq2}. For a given initial iterate $x_0$, the CG methods generates a sequence of iterates $\{x_k\}$ via
\begin{equation}\label{eq3}
	x_{k+1}=x_k+\lambda _kd_k, ~~~k=0,1,...,
\end{equation}
where $\lambda _k>0$ is a steplength, and $d_k$ is a search direction generated by
\begin{equation}\label{eq4}
	d_k=
	\begin{cases}
		-g(x_0), ~~~\text{if } k=0\\[10pt]
		-g(x_k)+\beta _kd_{k-1}, ~~~\text{if } k\geq 1.
	\end{cases}
\end{equation}
The classical Hestesnes-Stiefel (HS) conjugate gradient method \cite{hestenes1952} and Polak--Ribi{\`e}re--Polyak (PRP) \cite{polak1969,polyak1969} are popular methods for solving \eqref{eq2},
their parameters are respectively given by
\begin{equation}\label{eq5}
	\beta _k^{HS}=\frac{\langle g(x_k),~y_{k-1}\rangle}{\langle d_{k-1},~y_{k-1}\rangle}, \quad \beta _k^{PRP}=\frac{\langle g(x_k),~y_{k-1}\rangle}{\|g(x_{k-1})\|},
\end{equation}
where $y_{k-1}=g(x_k)-g(x_{k-1})$.

Various modifications of the HS and PRP methods for sloving \eqref{eq2} have been proposed in the literature (see, e.g.\cite{dong2015modified,amini2019modified,dong2017some,dong15,dong2020modified,babaie2016descent, gilbert1992,cheng2007,sun2015three,narushima2011,zhang2006descent,grippo1997globally,hager2005new,hager2006algorithm}). 
On the other hand, numerous methods on modified HS and PRP methods have been presented to address large-scale nonlinear systems of equations \eqref{eq1}, some of them can be found in the following references \cite{mohammad2021, awwal2020projection,cheng2009family,li2014derivative,tarzanagh2017,li2014,yu2010derivative,zhou2014inexact,zhou2015convergence}.

The multivariate spectral approach is a technique used to assemble (at each iteration step), a diagonal matrix approximation of the Jacobian of $F$. Based on our research findings, the study of the multivariate spectral approach for solving problem \eqref{eq1} or \eqref{eq2} has not received much attention in the literature. Han et al.\cite{han2008} presented the multivariate spectral method for unconstrained optimization. Yu et al. \cite{yu2011}, Liu and Li \cite{liu2015spectral} and Mohammad \cite{mohammad2021} extended this approach to solve large-scale nonlinear systems of monotone equations.  To the best of our knowledge, no existing literature on the multivariate spectral approach combined with the hybridization of the conjugate gradient method for solving general nonlinear systems of equations. Therefore, it is interesting to propose a multivariate spectral hybridization of the HS and PRP method for solving \eqref{eq1}.

The aim of this paper is to present a derivative-free method for solving large-scale general nonlinear systems of equations. This aim will be achieved through a positive hybridization of the HS and PRP methods and incorporating multivariate spectral approach. 

The remaining parts of this paper are as follows. In Section 2 we present the algorithm of the proposed method. In Section 3 we establish the global convergence of the proposed algorithm. In Section 4 we present the numerical experiments, and conclusions in Section 5. Unless otherwise stated, throughout this paper $\|\cdot\|$ stands for the Euclidean norm of vectors, and $\langle\,\cdot~,~\cdot\rangle$ is the inner product of vectors in $\mathbb{R}^n$. The merit function is given as
\begin{equation}\label{eq6}
	f(x):=\frac{1}{2}\|F(x)\|^2,
\end{equation}
which implies that  $f(x)=0$ if and only if $F(x)=0$.
\section{Algorithm}
In this section,  we present a positive hybridization of HS and PRP (hybrid+) method incorporating a diagonal approximation of the Jacobian of the function $F$ for solving \eqref{eq1}. The sequence of iterates is given by \eqref{eq3}
where the search direction $d_k$ in this case is generated by
\begin{equation}\label{eq7}
	d_k=
	\begin{cases}
		-F(x_0), ~~~\text{if } k=0\\[10pt]
		-B_kF(x_k)+\bar{\beta} _kd_{k-1}, ~~~\text{if } k\geq 1,
	\end{cases}
\end{equation}
where $B_k$ is a positive definite diagonal matrix given by 
\begin{equation}\label{eq8}
	B_k=\text{diag}\biggl(\frac{1}{b _k^1},\frac{1}{b _k^2},...,\frac{1}{b _k^n}\biggr),
\end{equation}
\begin{equation}\label{eq9}
	b_{k}^i=
	\begin{cases}
		\max \left\{\min \biggl\{ \frac{y _{k-1}^i}{s_{k-1}^i},~~ u\biggr\},~~\ell \right\},  & \text{ if }s_{k-1}^i\ne 0 \\
		1, &  \text{ if }s_{k-1}^i= 0,
	\end{cases}
\end{equation}
where $y_{k-1}=F(x_k)-F(x_{k-1}),~s_{k-1}=x_k-x_{k-1}, ~0<\ell <1 \leq u$ and $b_{k}^i~(i=1,2,\ldots , n)$  is the $i$th non-zero component of the diagonal matrix $B$ at the iterate $k$. The parameter $\bar{\beta} _k$, is defined as 

\begin{equation}\label{eq10}
	\bar{\beta} _k=\frac{\max\{0, ~~\langle F(x_k),~y_{k-1}\rangle\}}{\max\{\langle d_{k-1},~y_{k-1}\rangle,~~\|F(x_{k-1})\|^2\}}.
\end{equation}

However, it can be observed that the direction defined using \eqref{eq7}--\eqref{eq10} may not be of descent if the objective function $f$ is given by \eqref{eq6}. Based on this observation, we recall the nonmonotone derivative-free backtracking line search for computing the steplength. Let $\lambda_k$ satisfy the following condition

\begin{equation}\label{eq12}
	f(x_k+\lambda_{k}d_k)\leq C_k+\tau_k-\sigma\lambda_{k}^2\|d_k\|^2,
\end{equation}
where $\sigma \in (0,1)$ and $\{\tau_k\}$ is a suitable sequence of positive numbers satisfying  
\begin{equation}\label{eq13}
	\sum _{k=0}^{\infty}\tau_k\leq \tau <\infty,
\end{equation}
and $C_k$ is updated via
\begin{equation}
	C_{k+1} =\frac{\eta_kQ_{k}(C_{k}+\tau_{k})+f(x_{k+1})}{Q_{k+1}},
\end{equation}
with $Q_0=1,~C_0=f_0,~Q_{k+1}=\eta_{k}Q_{k}+1$.

It is remarkable that the application of similar line search condition \eqref{eq12} to the spectral residual method has been studied by Cheng and Li in \cite{chengli2009}. In this paper, we will use the derivative-free nonmonotone line search condition \eqref{eq12} to the multivariate spectral hybridization of HS and PRP method \eqref{eq7}--\eqref{eq10}.  Below are the detailed steps of the proposed algorithm.


\begin{algorithm*}
	\scriptsize
	\caption{{Diagonal hybrid HS+ and PRP+ method (\texttt{hybrid+})}}\label{algo_dprp}
	\vspace{\fill}
	\SetKwData{Left}{left}\SetKwData{This}{this}\SetKwData{Up}{up}
	\SetKwFunction{Union}{Union}\SetKwFunction{FindCompress}{FindCompress}
	\SetKwInOut{Input}{Input}\SetKwInOut{Output}{output}
	\SetKwInOut{Stepo}{Step 1}\SetKwInOut{Stept}{Step 2}
	\SetKwInOut{Steptw}{Step 3}
	\SetKwInOut{Stepth}{Step 4}\SetKwInOut{Steptt}{Step 5}
	\SetKwInOut{Stepf}{Step 6}\SetKwInOut{Stepfv}{Step 7}\SetKwInOut{Stepsx}{Step 8}
	\Input{Given $x_0 \in \mathbb{R}^n,~ \rho,~\sigma$,~ $\eta_{min}, ~\eta_{max}\in (0,1),~ \text{tol}>0$, ~$0<\ell \leq 1\le u$. Set $k=0$,~$B_0=I$,~ $C_0=f(x_0)$,~$Q_0=1$. Choose a sequence $\{\tau_k\}$ satisfying 
		\(
		\displaystyle\sum _{k=0}^{\infty}\tau_k\leq \tau <\infty.
		\) }
	\BlankLine
	\Stepo{
		Compute $\|F(x_k)\|$, \If {$\|F(x_k)\|\leq \text{tol}$,}{ stop.}
	}
	\Stept{
		\eIf {$k=0$,} {
			set $d_k:=-F(x_k)$\;
		}
		{	Compute $d_k$  using \eqref{eq7}--\eqref{eq10}.
		}}\label{stp2}
	\Steptw{ Set $\lambda_{k}=1$.}
	\Stepth{	\uIf {\begin{equation}\label{eqls1} f(x_k+\lambda_{k}d_k)\leq C_k+\tau_k-\sigma\lambda_{k}^2\|d_k\|^2,\end{equation}  } {
		set $x_{k+1}:=x_k+s_k$, ~$s_k=\lambda_{k}d_k$.
		}
		\uElseIf{ \begin{equation}\label{eqls2} f(x_k-\lambda_{k}d_k)\leq C_k+\tau_k-\sigma\lambda_{k}^2\|d_k\|^2,\end{equation}  } {
			set  $x_{k+1}:=x_k-s_k$, ~$s_k=\lambda_{k}d_k$.	
		}
	\Else{ Set $\lambda_k=\rho \lambda_{k}$ and go to Step $4$.} 
		 }
	 \Steptt{Choose $\eta_k\in [\eta_{min},~\eta_{max}]$ and compute
	 \begin{equation*}\label{eqnonmono}
	 	Q_{k+1}=\eta_kQ_k+1, \quad C_{k+1}=\frac{\eta_kQ_k(C_k+\tau_k)+f(x_{k+1})}{Q_{k+1}.}
	 \end{equation*}
 
}
	\Stepf{ Set $k=k+1$ and go to Step 1.}
\end{algorithm*}\DecMargin{0.2em}
\begin{remark}\label{rmk1}
		By the definition of $\{b_k^i\}$ in \eqref{eq8}, the diagonal matrix $B_k$ is positive definite and bounded. In fact, we have
		\[\ell \leq b_k^i\leq u ~~~\text{for all~~} i,~k. \] 
\end{remark}
\begin{remark}\label{rmk2}
	From equation \eqref{eq10} it is clear that $\bar{\beta_k}$ is a hybridization of the HS+ and PRP+ parameters. Therefore, the resulting algorithm is expected to inherit the good performance of the two positive parameters in practice. 
\end{remark}
The following Lemma is obtained from Lemma 2.2 in \cite{chengli2009}.
\begin{lemma}\label{lm0}
	The iterates generated by Algorithm \ref{algo_dprp} satisfy	
	\begin{equation}\label{eq11c}
	0\leq f(x_k)\leq C_k\leq P_k~~~ \text{for all}~~ k\geq 0,
	\end{equation}  where $P_k$ is define by
	\begin{equation}\label{eq12b}
		P_k=\frac{1}{k+1}\sum_{j=0}^{k}(f(x_j)+j\tau_{j-1}), ~~\tau_{-1}=0.
	\end{equation}
Moreover, 
\begin{equation}\label{eq13b}
	C_{k+1}\leq C_k+\tau_k.
\end{equation}
\end{lemma}

Since $\tau_k>0$, it is easy to see that for all sufficiently small $\lambda_k$, the line search conditions \eqref{eqls1} and \eqref{eqls2} hold.  Also from Lemma \ref{lm0} we have  $f(x_k)\leq C_k$. Thus, the step length  obtained from the  nonmonotone line search step in Algorithm \ref{algo_dprp} is well-defined.

\section{Convergence Analysis}
The following standard assumptions are essential in proving the global convergence of the proposed algorithm.
\begin{assumption}\label{ass1}
	The level set $\Theta_{0}=\{x\in \mathbb{R}^n : f(x)\leq f(x_0)+\tau\}$ is bounded and the solution set of problem \eqref{eq1} is non-empty, that is there exists $\bar{x}\in \mathbb{R}^n$ such that $F(\bar{x})=0$.
\end{assumption}
\begin{assumption}\label{ass2}
	The function $F:\mathbb{R}^n\to \mathbb{R}^n$ is continuously differentiable and its Jacobian is Lipschitz continuous on $\tilde{\Theta} \supseteq \Theta_0$, i.e., there exists a constant $\mathcal{L}>0$ such that
	\begin{equation}\label{eq14}
		\|J(x)-J(y)\|\leq\mathcal{L}\|x-y\|, ~~~\forall x,~y\in \tilde{\Theta}.
	\end{equation}
\end{assumption}
Assumption \ref{ass2} implies that there exist nonnegative constants $\omega_1,~~\omega_2$ and $\mathcal{L}_1$
such that
\begin{equation}\label{eq15}
	\|F(x)\|\leq \omega_1, ~~~~\|J(x)\|\leq \omega_2 ~~~~\text{for all } x\in \tilde{\Theta}.
\end{equation}
\begin{equation}\label{eq15b}
	\|F(x)-F(y)\|\leq \mathcal{L}_1\|x-y\|, ~~~~\|\nabla f(x)-\nabla f(y)\|\leq \mathcal{L}_1\|x-y\| ~~~~\text{for all } x,~y \in \tilde{\Theta}.
\end{equation}
The proof of the following Lemma is similar to Lemma 3.1 in \cite{chengli2009}.
\begin{lemma}\label{lm1}
	Suppose Assumption \ref{ass1} holds. Then the sequence $\{x_k\}$ generated by Algorithm \ref{algo_dprp} is contained in the level set $\Theta_0$.
\end{lemma}
\begin{lemma}\label{lm3}
	Let $\{x_k\}$ be the sequence generated by Algorithm \ref{algo_dprp}, then we have 
	\[\displaystyle\lim_{k\to \infty}\lambda_{k}\|d_k\|=0.\]
\end{lemma}
\begin{proof}
 From Step $4$ of Algorithm \ref{algo_dprp}  for any $k\geq 0$ we have
	\begin{align}\label{eq16}
	f(x_{k+1})	 & \leq 
		C_k +\tau_k-\sigma \lambda_{k}^2\|d_k\|^2.
	\end{align}
	Combining \eqref{eqnonmono} and  \eqref{eq16}, we have
	\begin{align*}
		C_{k+1} & =\frac{\eta_kQ_k(C_k+\tau_k)+f(x_{k+1})}{Q_{k+1}}\\
		&\leq \frac{\eta_kQ_k(C_k+\tau_k)+C_k +\tau_k-\sigma \lambda_{k}^2\|d_k\|^2}{Q_{k+1}}\\
		&= C_k+\tau_k-\frac{\sigma \lambda_{k}^2\|d_k\|^2}{Q_{k+1}}.
		\end{align*}
	Since  $\displaystyle\sum_{k=0}^{\infty}\tau_k<~\infty$  and $\sigma$ is a positive constant, it follows that 
	\begin{equation}\label{eq17}
		\sum_{k=0}^{\infty} \frac{\lambda_{k}^2\|d_k\|^2}{Q_{k+1}}<~\infty. 
	\end{equation}
From \eqref{eqnonmono} and the fact that $\eta_{max}\in (0,1)$, $Q_0=1$, it follows that

	\[Q_{k+1} =1+\eta_kQ_k=1+\sum_{i=0}^{k}\prod_{j=0}^{i}\eta_{k-j}\leq 1+\sum_{i=0}^{k}\eta_{max}^{i+1}\leq \sum_{i=0}^{\infty}\eta_{max}^{i}=\frac{1}{1-\eta_{max}}.\]
This and \eqref{eq17} gives
\begin{equation}\label{eq18}
		\lim_{k\to \infty}\lambda_{k}\|d_k\|=0.
	\end{equation}
\end{proof}
\begin{theorem}\label{thm1}
	Let Assumption \ref{ass2} holds, if $\{x_k\}$ is the sequence generated by Algorithm \ref{algo_dprp}. Then either
	\begin{equation}\label{eq19}
		\liminf_{k\to \infty}\|F(x_k)\|=0,
	\end{equation}
or 
\begin{equation}\label{eq20}
	\liminf_{k\to \infty}\langle J(x_k)^TF(x_k),~B_kF(x_k)\rangle= 0.\end{equation}
	
\end{theorem}
\begin{proof}
	Suppose that \eqref{eq19} is not true. That is, there exists $\xi>0$ such that 
	\begin{equation}\label{eq21} 
		\|F(x_k)\|\geq \xi~~~\text{for all}~~k\in \mathbb{N}\cup \{0\}.
	\end{equation}
By the definition of the search direction $d_k$ in \eqref{eq7}, if $k=0$ then by \eqref{eq15} we get
	\begin{equation}\label{eq21b}
		\|d_k\|=\|F(x_k)\|\leq \omega_1.
\end{equation}
Otherwise, for $k\geq 1$. If $\bar{\beta_k}= 0$, then  
\begin{align*}
	\|d_k\|^2&=\|B_kF(x_k)\|^2\\
	&= \sum_{i=1}^nF(x_k)^i(b_k^i)^{-2}F(x_k)^i\\
	&\leq \frac{1}{\ell^2}\sum_{i=1}^{n}(F(x_k)^i)^2\\
	&=\frac{1}{\ell^2}\|F(x_k)\|^2,
\end{align*} 
where $F(x_k)^i$ is the $i$th component of the vector $F$ at $x_k$.
This and \eqref{eq15}, gives
\begin{equation}\label{eq22}
	\|d_k\|\leq \frac{\omega_1}{\ell}.
\end{equation}
If  $\bar{\beta_k}>0$, then from \eqref{eq10}, Cauchy-Schwartz inequality, \eqref{eq15}, \eqref{eq15b}, \eqref{eq18} and \eqref{eq21}, we have
\begin{equation}\label{eq22b}
	|\bar{\beta_k}|\leq \frac{\|F(x_k)\|\|y_{k-1}\|}{\|F(x_{k-1})\|^2}\leq \frac{\mathcal{L}_1\omega_1\lambda_{k-1}\|d_{k-1}\|}{\xi^2}\to 0.
\end{equation}
This implies that there exists $\mu\in (0,~1)$ such that for sufficiently large $k$, 
\begin{equation}\label{eq23}
	|\bar{\beta_k}|\leq \mu.
\end{equation}
Suppose without loss of generality that \eqref{eq23} holds for all $k$, then from \eqref{eq7} and triangular inequality it follows that
\begin{align}
	\|d_k\|&\leq \|B_kF(x_k)\|+|\bar{\beta_k}|\|d_{k-1}\|\nonumber\\
	&\leq \frac{\omega_1}{\ell}+\mu\|d_{k-1}\|\nonumber\\
	&\leq \frac{\omega_1}{\ell}(1+\mu+\mu^2+\ldots+\mu^{k-1})+\mu^k\|d_0\|\nonumber\\
	&= \frac{\omega_1}{\ell}(1+\mu+\mu^2+\ldots+\mu^{k-1})+\mu^k\|F_0\|\nonumber\\
	&\leq \frac{\omega_1}{\ell(1-\mu)}+\omega_1\label{eq24}.
\end{align}
Inequalities \eqref{eq21b}, \eqref{eq22} and \eqref{eq24} give
\begin{equation}\label{eq25}
	\|d_k\|\leq N,
\end{equation}
where $N=\max\left\{\omega_1,~\frac{\omega_1}{\ell},~\frac{\omega_1}{\ell(1-\mu)}+\omega_1\right\}.$\\[5mm]
Next, we consider the following two cases:\\[2mm]
(i.) If $\displaystyle\liminf_{k\to \infty}\|d_k\|=0,$ then by definition of the search direction $d_k$ and the boundedness of the diagonal matrix $B_k$ we have
	\begin{gather*}
		\frac{1}{u}\|F(x_k)\|\leq \|B_kF(x_k)\|= \|d_k\|+|\bar{\beta_k}|\|d_{k-1}\|,
	\end{gather*}
which together with \eqref{eq22b} and the assumption that $\displaystyle\liminf_{k\to \infty}\|d_k\|=0$ gives
\begin{align*}
	\frac{1}{u}\liminf_{k\to \infty}\|F(x_k)\|\leq \liminf_{k\to \infty}\|d_k\|+\liminf_{k\to \infty}|\bar{\beta_k}|\|d_{k-1}\|=0.
\end{align*}
It follows that $\displaystyle\liminf_{k\to \infty}\|F(x_k)\|=0,$ which contradicts \eqref{eq21}.\\
(ii.) If $\displaystyle\liminf_{k\to \infty}\|d_k\|>0,$ then by \eqref{eq18} we have
\begin{equation}\label{eq26}
	\liminf_{k\to \infty}\lambda_{k}=0.
\end{equation} 
From Step $4$ of Algorithm \ref{algo_dprp}, for sufficiently large $k$,~ $\frac{\lambda _k}{\rho}$ does not satisfy \eqref{eqls1} and \eqref{eqls2}, that is 
\begin{equation}\label{eq27a}
	f(x_k+\frac{\lambda_k} {\rho} d_k) >C_k+  \tau_k-\sigma\left(\frac{\lambda_k} {\rho}\right)^2\|d_k\|^2,
	\end{equation} 
and
\begin{equation}\label{eq27b}
	f(x_k-\frac{\lambda_k} {\rho} d_k) >C_k+  \tau_k-\sigma\left(\frac{\lambda_k} {\rho}\right)^2\|d_k\|^2. 
\end{equation}
Subtituting \eqref{eq11c} in \eqref{eq27a}, we obtain
\begin{equation}\label{eq27c}
	\frac{f\left(x_k+\frac{\lambda_k}{\rho} d_k\right)-f(x_k)}{\lambda_k}>-\sigma\rho^{-2}\lambda_k\|d_k\|^2.\end{equation}
By mean value theorem, there exists $\varepsilon\in (0,1)$ such that 
\begin{equation}\label{eq27d}
		\frac{f\left(x_k+\frac{\lambda_k}{\rho} d_k\right)-f(x_k)}{\lambda_k}=\rho^{-1}\left\langle J\left(x_k+\varepsilon\frac{\lambda_k}{\rho}d_k\right)^TF\left(x_k+\varepsilon\frac{\lambda_k}{\rho}d_k\right),~d_k\right\rangle.
\end{equation}
Applying Cauchy-Schwartz inequality and \eqref{eq15b} on the right hand side of \eqref{eq27d} yields
	\begin{align}\label{eq28}
		\rho^{-1}\left\langle J\left(x_k+\varepsilon\frac{\lambda_k}{\rho}d_k\right)^TF\left(x_k+\varepsilon\frac{\lambda_k}{\rho}d_k\right),~d_k\right\rangle & =	\rho^{-1}\left\langle J(x_k)^TF(x_k),~d_k\right\rangle\nonumber\\ & +	\rho^{-1}\left\langle J\left(x_k+\varepsilon\frac{\lambda_k}{\rho}d_k\right)^TF\left(x_k+\varepsilon\frac{\lambda_k}{\rho}d_k\right)-J(x_k)^TF(x_k),~d_k\right\rangle\nonumber\\
		& \leq 	\rho^{-1}\left\langle J(x_k)^TF(x_k),~d_k\right\rangle\nonumber\\
		&+ 	\rho^{-1}\left\|J\left(x_k+\varepsilon\frac{\lambda_k}{\rho}d_k\right)^TF\left(x_k+\varepsilon\frac{\lambda_k}{\rho}d_k\right)-J(x_k)^TF(x_k)\right\|\|d_k\|\nonumber\\
		&\leq 	\rho^{-1}\langle J(x_k)^TF(x_k),~d_k\rangle +\mathcal{L}_1\varepsilon\frac{\lambda_k}{\rho^2}\|d_k\|^2.
\end{align}
Combining \eqref{eq27c} and \eqref{eq28} and using \eqref{eq7} together with \eqref{eq25}, we have
\begin{equation*}
	-(\sigma+\mathcal{L}_1\varepsilon)\frac{\lambda_k}{\rho}\|d_k\|^2<\langle J(x_k)^TF(x_k),~d_k\rangle. \label{eq28a}
\end{equation*}
This imply \begin{equation}\label{eq28b}
	-(\sigma+\mathcal{L}_1\varepsilon)\frac{\lambda_k}{\rho}N^2<-\langle J(x_k)^TF(x_k),~B_kF(x_k)\rangle+\bar{\beta_k}\langle J(x_k)^TF(x_k),~d_{k-1}\rangle.\end{equation}
It follows from \eqref{eq28b} that
\begin{equation}\label{eq29}
\liminf_{k\to \infty}\langle J(x_k)^TF(x_k),~B_kF(x_k)\rangle\leq 0.
\end{equation}
Using \eqref{eq27b} and repeating the above process in a similar manner, yields
\begin{equation}\label{eq30}
	\liminf_{k\to \infty}\langle J(x_k)^TF(x_k),~B_kF(x_k)\rangle\geq 0.
\end{equation}
 The inqualities \eqref{eq29} and \eqref{eq30} imply equation \eqref{eq20}.
\end{proof}

\section{Numerical Results}\label{sec4}
In this section, we compare the computational performance of the proposed method \texttt{hybrid+} with the following existing methods for solving general nonlinear equations.
\begin{itemize}
	\item \texttt{dfriml}: A new derivative-free conjugate gradient method for large-scale nonlinear equations proposed in \cite{fang2017}. 
	\item \texttt{dfprp}: A derivative-free PRP method for solving large-scale nonlinear systems of equations and its global convergence proposed in \cite{li2014derivative}. 
\end{itemize}
The parameters used for the implementation of the \texttt{dfriml} and \texttt{dfprp}  methods are set just as in \cite{fang2017} and \cite{li2014derivative} respectively.

The parameters used for the implemetation of the proposed  \texttt{hybrid+} method are specified as follows: $\rho=0.5,~\sigma =10^{-4}$, $\eta_{\text{min}}=0.1, ~\eta_{\text{max}}=0.85,$~$\ell =10^{-10}, ~u=10^{10}$, and $\tau_k=\frac{1}{2^{k}}.$ The parameter $\eta_k$ is chosen as
\[\eta_k=0.75e^{\min\{\hat{\omega},~(\frac{k}{75})^2\}}+0.1, ~~0<\hat{\omega}<0.18.\] 

We tested three methods on the following ten test problems, where $F=(F_1,F_2, ... ,F_n)^T.$\\[2mm]

\noindent \emph{Problem 1}: Modified exponential function \cite{lacruz2004}
\begin{align*}
	F_1(x)&=e^{x_1}-1\\
	F_i(x) &= e^{x_i}+x_{i}-1, ~~ i=2,3\ldots,n-1.
\end{align*}

\noindent \emph{Problem 2}:  Logarithmic function \cite{lacruz2004}
\begin{equation*}
	F_i(x_i)=\log (x_i+1)-\frac{x_i}{n}, ~~ i=1,2,\ldots,n.
\end{equation*}

\noindent \emph{Problem 3}: Strictly convex function I \cite{lacruz2004}
\begin{equation*}
	F_i(x)=e^{x_i}-1, ~~ i=1,2,\ldots,n.
\end{equation*}
\noindent \emph{Problem 4}: Modified strictly convex function II \cite{lacruz2004}
\begin{equation*}
	F_i(x)=\left(\dfrac{i}{n+1}\right)e^{x_i}-1, ~~ i=1,2,\ldots,n.
\end{equation*}

\noindent \emph{Problem 5}: Tridiagonal exponential function \cite{lacruz2004}
\begin{align*}
	F_1(x)& =x_1-e^{\cos (h(x_1+x_2))}\\
	F_i(x)&=x_i-e^{\cos (h(x_{i-1}+x_i+x_{i+1}))}, ~~ i=2,\ldots,n-1,\\
	F_n(x)&=x_n-e^{\cos (h(x_{n-1}+x_n))},\\
	h&=\frac{1}{n+1}.
\end{align*}
\noindent \emph{Problem 6}: Gradient of Engval function \cite{li1999}
\begin{align*}
	F_1(x)&=x_1(x_1^2+x_2^2)-1\\
	F_i(x)&=x_i(x_{i-1}^2+2x_i^2+x_{i+1}^2)-1,~~2\leq i \leq n-1\\
	F_n(x)&=x_n(x_{n-1}^2+x_n^2).
\end{align*}
\noindent \emph{Problem 7}: Chandrasekhar H-equation \cite{kelly1980} 
\begin{equation*}
	F_i(x) = x_i - \left(1-\dfrac{c}{2n}\sum_{j=1}^{n}\dfrac{\delta_ix_j}{\delta_i+\delta_j}\right)^{-1},\quad
	c=0.9,~~ \delta_i=\dfrac{i-0.5}{n},~~i=1,2,\ldots,n.
\end{equation*}
\noindent \emph{Problem 8}: Modified Problem $3.34$ in  \cite{lukvsan2018}
\begin{align*}
	F_i(x)&=x_i-\frac{x_{i+1}^3}{100},~~1\leq i \leq n-1\\
	F_n(x)&=x_n-\frac{x_{n}^3}{100}.
\end{align*}
\noindent \emph{Problem 9}: Nonsmooth function 1 \cite{yu09} 
\begin{align*}
	F_{i}(x)&=x_{i}-\sin(|x_i-1|),  ~ \text{for }i=1,2,\ldots,n.
	\end{align*}
\noindent \emph{Problem 10}: Nonsmooth function 2 \cite{li2011}
\begin{align*}
	F_{i}(x)&=2x_{i}-\sin(|x_i|),  ~ \text{for }i=1,2,\ldots,n.
\end{align*}

We used five different dimensions $n=1000,~5000,~10000,~50000,~100000$ with  ten different initial points as follows:\\ $x_{0}^1=(1,\ldots,1)^T$, ~$x_{0}^2=(0.1,\ldots ,0.1)^T$, ~$x_{0}^3=\bigl(\frac{1}{2},\ldots ,\frac{1}{2^n}\bigr)^T$,  ~$x_{0}^4=\bigl(1-\frac{1}{n}, 1-\frac{2}{n},\ldots,0)^T$,
$x_{0}^5=\bigl(0, \frac{1}{n},\frac{2}{n},\ldots ,\frac{n-1}{n}\bigr)^T$, ~$x_{0}^6=\bigl(1,\frac{1}{2},\ldots ,\frac{1}{n}\bigr)^T$, ~$x_{0}^7=\bigl(\frac{n-1}{n}, \frac{n-2}{n},\ldots ,0 \bigr)^T$, ~$x_{0}^8=\bigl(\frac{1}{n}, \frac{2}{n},\ldots ,1\bigr)^T$, $x_{0}^9=(10,10,\ldots ,10)^T$  and $x_{0}^{10}=\text{rand}(0,1)$. Here, $\text{rand}(0,1)$ means the initial point is chosen randomly from the interval $(0,1)$. 
All the three methods are coded in MATLAB R2018a and run on a PC with an intel COREi7 processor, 8GB of RAM and CPU 1.8GHz $\times$ 8. 
All implementation processes are is stopped when $\|F(x_k)\|\le 10^{-6}$ or number of iterations exceeds $1000.$

In Tables \ref{tab1}--\ref{tab10} (see Appendix A) we present  results on the following information: the number of iterations (\#iter) needed to converge to an approximate solution, the number of function evaluations (\#fval), the CPU time in seconds (time) and the norm of the function at the approximate solution (Fnorm).

\begin{figure}[htp]
	\begin{center}
		\includegraphics[angle=0,width=12cm]{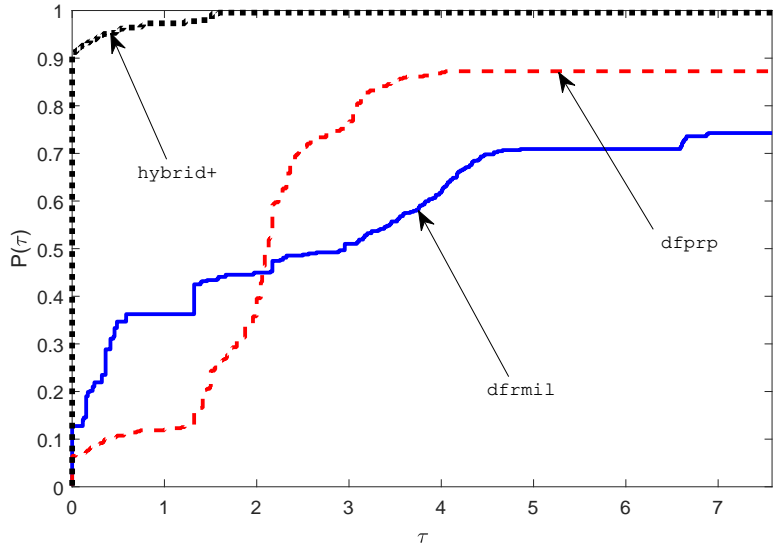}
	\end{center}
	\vspace{-0.5cm}
	{{\caption{Dolan and Mor$\acute{e}$ performance profile with respect to number of iterations }\label{fig1}}}
\end{figure}

\begin{figure}[htp]
	\begin{center}
		\includegraphics[angle=0,width=12cm]{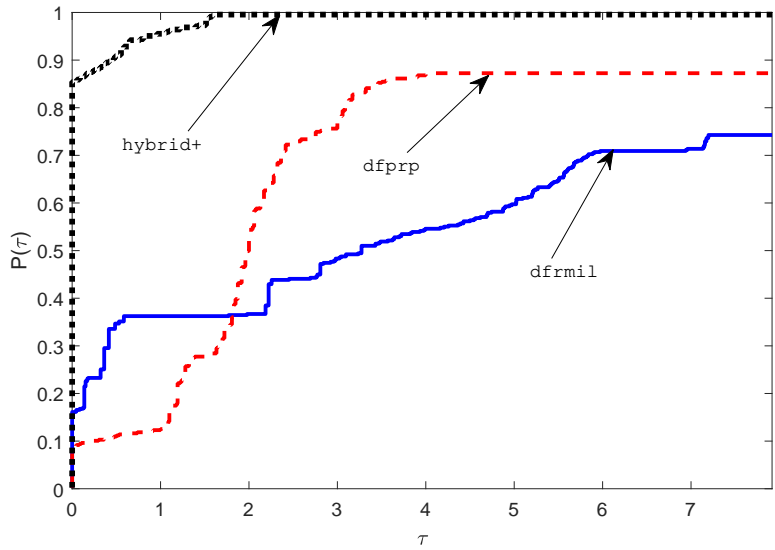}
	\end{center}
	\vspace{-0.5cm}
	{{\caption{Dolan and Mor$\acute{e}$ performance profile with respect to number of function evaluations}\label{fig2}}}
\end{figure}

\begin{figure}[htp]
	\begin{center}
		\includegraphics[angle=0,width=12cm]{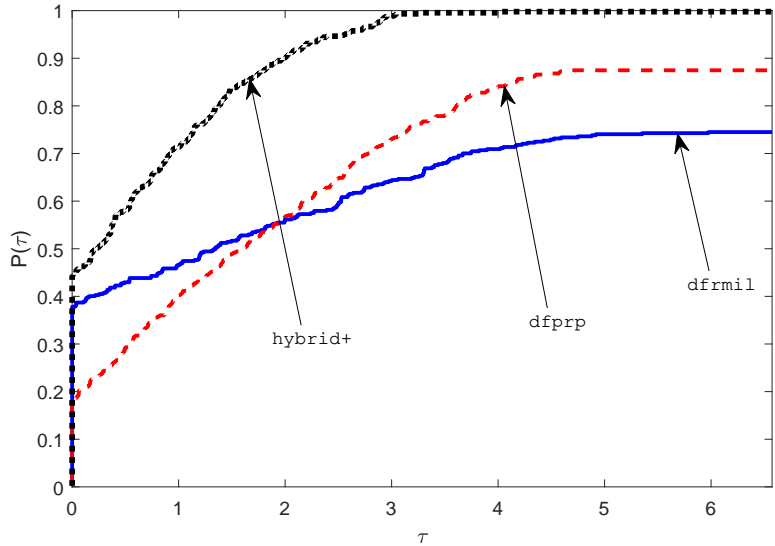}
	\end{center}
	\vspace{-0.5cm}
	{{\caption{Dolan and Mor$\acute{e}$ performance profile with respect to time}\label{fig3}}}
\end{figure}

Figures \ref{fig1}--\ref{fig3} shows the performance of the tested methods relative to \#iter, \#fval and time. Efficiency comparisons were made using the performance profile introduced by Dolan and Mor$\acute{e}$ \cite{dolan2002}.

It can be observed from Figures \ref{fig1}--\ref{fig3} that \texttt{hybrid+} is the best solver since it becomes the most efficient and robust method based on the performed numerical experiments. The performance of \texttt{hybrid+} in terms of number of iterations metric (Fig. \ref{fig1}) is commendable since it solves and wins more than  $90\%$ of the problems with few number of iterations. Results from Fig. \ref{fig2} and Fig. \ref{fig3}  show that \texttt{hybrid+} is very efficient since it solves more than $80\%$ of the problems with few number of function evaluations, and  has less CPU time compared to the \texttt{dfrmil} and \texttt{dfprp} methods.

\newpage
\section{Conclusion}
We proposed a multivariate spectral gradient method based on the hybridization of HS and PRP methods for solving nonlinear systems of equations. The interesting feature of the method is that: it does not require derivative as well as matrix storage. Thus, it is simple and easy to implement.  The global convergence of the sequence generated by the method is obtained using a derivative-free nonmonotone line search technique. The numerical comparison of the proposed method with similar methods in the literature indicates that it is promising.


\bibliographystyle{plain}
\bibliography{multivariatehybrid} 
\appendix
\section*{Appendix A:  Table of Numerical Experiments}\medskip
Below is the details of the numerical experiments conducted in Section \ref{sec4}.
\begin{table}[htbp]
	\centering
	\caption{Test results for the \texttt{hybrid+}, \texttt{dfrmil} and \texttt{dfprp} methods on problem 1}
	\resizebox{\columnwidth}{!}{
\begin{tabular}{|c|r|c|c|c|c|c|c|c|c|c|c|c|c|}
			\hline
			\multirow{2}[4]{*}{Dimension} & \multicolumn{1}{c|}{\multirow{2}[4]{*}{$x_0$}} & \multicolumn{4}{c|}{dfrmil}   & \multicolumn{4}{c|}{dfprp}    & \multicolumn{4}{c|}{hybrid+} \bigstrut\\
			\cline{3-14}          & \multicolumn{1}{c|}{} & \#iter & \#fval & time  & Fnorm & \#iter & \#fval & time  & Fnorm & \#iter & \#fval & time  & Fnorm \bigstrut\\
			\hline
			\multirow{10}[20]{*}{1000} & $x_0^1$ & 686   & 1204  & 0.32792 & 9.60E-07 & F     & F     & F     & F     & 9     & 11    & 0.068379 & 4.78E-07 \bigstrut\\
			\cline{2-14}          & $x_0^2$ & 445   & 775   & 0.26144 & 7.99E-07 & 22    & 22    & 0.027148 & 7.60E-07 & 7     & 8     & 0.029527 & 1.86E-08 \bigstrut\\
			\cline{2-14}          & $x_0^3$ & F     & F     & F     & F     & 25    & 25    & 0.020625 & 6.49E-07 & 12    & 13    & 0.019615 & 6.78E-10 \bigstrut\\
			\cline{2-14}          & $x_0^4$ & F     & F     & F     & F     & F     & F     & F     & F     & 8     & 10    & 0.019977 & 5.58E-09 \bigstrut\\
			\cline{2-14}          & $x_0^5$ & F     & F     & F     & F     & F     & F     & F     & F     & 7     & 9     & 0.024601 & 3.76E-07 \bigstrut\\
			\cline{2-14}          & $x_0^6$ & F     & F     & F     & F     & 25    & 25    & 0.010021 & 8.96E-07 & 10    & 14    & 0.026651 & 1.88E-08 \bigstrut\\
\hline
		\end{tabular}%
	}
	\label{tab1}%
\end{table}%

\begin{table}[htbp]
	\centering
	\caption{Test results for the \texttt{hybrid+}, \texttt{dfrmil} and \texttt{dfprp} methods on problem 2}
	\resizebox{\columnwidth}{!}{
		\begin{tabular}{|c|r|c|c|c|c|c|c|c|c|c|c|c|c|}
			\hline
			\multirow{2}[4]{*}{Dimension} & \multicolumn{1}{c|}{\multirow{2}[4]{*}{$x_0$}} & \multicolumn{4}{c|}{dfrmil}   & \multicolumn{4}{c|}{dfprp}    & \multicolumn{4}{c|}{hybrid+} \bigstrut\\
			\cline{3-14}          & \multicolumn{1}{c|}{} & \#iter & \#fval & time  & Fnorm & \#iter & \#fval & time  & Fnorm & \#iter & \#fval & time  & Fnorm \bigstrut\\
			\hline
			\multirow{10}[20]{*}{1000} & $x_0^1$ & F     & F     & F     & F     & 32    & 32    & 0.025396 & 9.53E-07 & 6     & 6     & 0.0139 & 2.58E-07 \bigstrut\\
			\cline{2-14}          & $x_0^2$ & F     & F     & F     & F     & 28    & 28    & 0.03116 & 8.99E-07 & 4     & 4     & 0.008667 & 2.11E-09 \bigstrut\\
			\cline{2-14}          & $x_0^3$ & 7     & 7     & 0.007411 & 6.85E-07 & 25    & 25    & 0.015634 & 8.19E-07 & 5     & 5     & 0.018212 & 3.01E-07 \bigstrut\\
			\cline{2-14}          & $x_0^4$ & 9     & 9     & 0.00324 & 5.45E-07 & 32    & 32    & 0.025705 & 6.06E-07 & 6     & 6     & 0.010801 & 6.06E-07 \bigstrut\\
			\cline{2-14}          & $x_0^5$ & 9     & 9     & 0.002939 & 5.45E-07 & 32    & 32    & 0.015773 & 6.06E-07 & 6     & 6     & 0.007913 & 6.06E-07 \bigstrut\\
			\cline{2-14}          & $x_0^6$ & 8     & 8     & 0.02  & 1.84E-07 & 27    & 27    & 0.012899 & 6.99E-07 & 6     & 6     & 0.008796 & 1.27E-09 \bigstrut\\
\hline
		\end{tabular}%
	}
	\label{tab2}%
\end{table}%

\begin{table}[htbp]
	\centering
	\caption{Test results for the \texttt{hybrid+}, \texttt{dfrmil} and \texttt{dfprp} methods on problem 3}
	\resizebox{\columnwidth}{!}{
		\begin{tabular}{|c|r|c|c|c|c|c|c|c|c|c|c|c|c|}
			\hline
			\multirow{2}[4]{*}{Dimension} & \multicolumn{1}{c|}{\multirow{2}[4]{*}{$x_0$}} & \multicolumn{4}{c|}{dfrmil}   & \multicolumn{4}{c|}{dfprp}    & \multicolumn{4}{c|}{hybrid+} \bigstrut\\
			\cline{3-14}          & \multicolumn{1}{c|}{} & \#iter & \#fval & time  & Fnorm & \#iter & \#fval & time  & Fnorm & \#iter & \#fval & time  & Fnorm \bigstrut\\
			\hline
			\multirow{10}[20]{*}{1000} & $x_0^1$ & F     & F     & F     & F     & 32    & 32    & 0.023338 & 7.49E-07 & 7     & 8     & 0.010324 & 4.51E-07 \bigstrut\\
			\cline{2-14}          & $x_0^2$ & F     & F     & F     & F     & 28    & 28    & 0.015069 & 8.84E-07 & 5     & 5     & 0.006591 & 5.46E-11 \bigstrut\\
			\cline{2-14}          & $x_0^3$ & 8     & 8     & 0.005716 & 7.24E-09 & 25    & 25    & 0.009789 & 7.71E-07 & 6     & 6     & 0.007827 & 1.29E-07 \bigstrut\\
			\cline{2-14}          & $x_0^4$ & 9     & 10    & 0.004666 & 6.93E-07 & 31    & 31    & 0.009446 & 8.82E-07 & 9     & 10    & 0.010306 & 1.65E-13 \bigstrut\\
			\cline{2-14}          & $x_0^5$ & 9     & 10    & 0.015547 & 6.93E-07 & 31    & 31    & 0.019485 & 8.82E-07 & 9     & 10    & 0.01305 & 1.67E-13 \bigstrut\\
			\cline{2-14}          & $x_0^6$ & 10    & 10    & 0.005791 & 9.36E-09 & 26    & 26    & 0.018174 & 9.55E-07 & 8     & 8     & 0.010774 & 5.26E-08 \bigstrut\\
\hline
		\end{tabular}%
	}
	\label{tab3}%
\end{table}%

\begin{table}[htbp]
	\centering
	\caption{Test results for the \texttt{hybrid+}, \texttt{dfrmil} and \texttt{dfprp} methods on problem 4}
	\resizebox{\columnwidth}{!}{
		\begin{tabular}{|c|r|c|c|c|c|c|c|c|c|c|c|c|c|}
			\hline
			\multirow{2}[4]{*}{Dimension} & \multicolumn{1}{c|}{\multirow{2}[4]{*}{$x_0$}} & \multicolumn{4}{c|}{dfrmil}   & \multicolumn{4}{c|}{dfprp}    & \multicolumn{4}{c|}{hybrid+} \bigstrut\\
			\cline{3-14}          & \multicolumn{1}{c|}{} & \#iter & \#fval & time  & Fnorm & \#iter & \#fval & time  & Fnorm & \#iter & \#fval & time  & Fnorm \bigstrut\\
			\hline
			\multirow{10}[20]{*}{1000} & $x_0^1$ & 15    & 15    & 0.007262 & 2.21E-08 & 31    & 31    & 0.007125 & 9.46E-07 & 12    & 22    & 0.024729 & 7.84E-08 \bigstrut\\
			\cline{2-14}          & $x_0^2$ & 10    & 10    & 0.004264 & 1.22E-09 & 31    & 31    & 0.022866 & 8.43E-07 & 8     & 9     & 0.007088 & 1.85E-09 \bigstrut\\
			\cline{2-14}          & $x_0^3$ & 9     & 9     & 0.003674 & 5.31E-08 & 31    & 31    & 0.020309 & 7.17E-07 & 7     & 8     & 0.009625 & 9.83E-10 \bigstrut\\
			\cline{2-14}          & $x_0^4$ & 11    & 12    & 0.006169 & 1.15E-07 & 32    & 32    & 0.011652 & 7.16E-07 & 8     & 9     & 0.012406 & 8.27E-11 \bigstrut\\
			\cline{2-14}          & $x_0^5$ & 15    & 15    & 0.012637 & 3.45E-10 & 31    & 31    & 0.022339 & 9.56E-07 & 11    & 21    & 0.014863 & 3.59E-07 \bigstrut\\
			\cline{2-14}          & $x_0^6$ & 9     & 9     & 0.003975 & 5.55E-08 & 31    & 31    & 0.020498 & 7.23E-07 & 9     & 10    & 0.01395 & 4.02E-10 \bigstrut\\
\hline
		\end{tabular}%
	}
	\label{tab4}%
\end{table}%

\begin{table}[htbp]
	\centering
	\caption{Test results for the \texttt{hybrid+}, \texttt{dfrmil} and \texttt{dfprp} methods on problem 5}
	\resizebox{\columnwidth}{!}{
		\begin{tabular}{|c|r|c|c|c|c|c|c|c|c|c|c|c|c|}
			\hline
			\multirow{2}[4]{*}{Dimension} & \multicolumn{1}{c|}{\multirow{2}[4]{*}{$x_0$}} & \multicolumn{4}{c|}{dfrmil}   & \multicolumn{4}{c|}{dfprp}    & \multicolumn{4}{c|}{hybrid+} \bigstrut\\
			\cline{3-14}          & \multicolumn{1}{c|}{} & \#iter & \#fval & time  & Fnorm & \#iter & \#fval & time  & Fnorm & \#iter & \#fval & time  & Fnorm \bigstrut\\
			\hline
			\multirow{10}[20]{*}{1000} & $x_0^1$ & 4     & 4     & 0.006695 & 1.63E-07 & 34    & 34    & 0.032489 & 6.15E-07 & 4     & 4     & 0.017199 & 1.27E-07 \bigstrut\\
			\cline{2-14}          & $x_0^2$ & 4     & 4     & 0.002269 & 1.89E-07 & 34    & 34    & 0.020395 & 9.37E-07 & 4     & 4     & 0.008289 & 1.47E-07 \bigstrut\\
			\cline{2-14}          & $x_0^3$ & 4     & 4     & 0.003074 & 1.89E-07 & 34    & 34    & 0.015061 & 9.72E-07 & 4     & 4     & 0.007991 & 1.52E-07 \bigstrut\\
			\cline{2-14}          & $x_0^4$ & 4     & 4     & 0.003293 & 1.80E-07 & 34    & 34    & 0.032448 & 8.01E-07 & 3     & 3     & 0.006362 & 1.45E-08 \bigstrut\\
			\cline{2-14}          & $x_0^5$ & 4     & 4     & 0.006365 & 1.80E-07 & 34    & 34    & 0.019338 & 8.01E-07 & 3     & 3     & 0.005618 & 1.45E-08 \bigstrut\\
			\cline{2-14}          & $x_0^6$ & 4     & 4     & 0.002818 & 1.89E-07 & 34    & 34    & 0.032067 & 9.70E-07 & 3     & 3     & 0.004033 & 9.39E-07 \bigstrut\\
\hline
		\end{tabular}%
	}
	\label{tab5}%
\end{table}%

\begin{table}[htbp]
	\centering
	\caption{Test results for the \texttt{hybrid+}, \texttt{dfrmil} and \texttt{dfprp} methods on problem 6}
	\resizebox{\columnwidth}{!}{
		\begin{tabular}{|c|r|c|c|c|c|c|c|c|c|c|c|c|c|}
			\hline
			\multirow{2}[4]{*}{Dimension} & \multicolumn{1}{c|}{\multirow{2}[4]{*}{$x_0$}} & \multicolumn{4}{c|}{dfrmil}   & \multicolumn{4}{c|}{dfprp}    & \multicolumn{4}{c|}{hybrid+} \bigstrut\\
			\cline{3-14}          & \multicolumn{1}{c|}{} & \#iter & \#fval & time  & Fnorm & \#iter & \#fval & time  & Fnorm & \#iter & \#fval & time  & Fnorm \bigstrut\\
			\hline
			\multirow{10}[20]{*}{1000} & $x_0^1$ & 662   & 2329  & 0.41455 & 9.01E-07 & F     & F     & F     & F     & 41    & 51    & 0.061731 & 7.42E-07 \bigstrut\\
			\cline{2-14}          & $x_0^2$ & 762   & 2605  & 0.51284 & 6.64E-07 & F     & F     & F     & F     & 40    & 51    & 0.048521 & 4.56E-07 \bigstrut\\
			\cline{2-14}          & $x_0^3$ & 667   & 2512  & 0.18062 & 9.64E-07 & F     & F     & F     & F     & 37    & 53    & 0.051431 & 8.01E-07 \bigstrut\\
			\cline{2-14}          & $x_0^4$ & 473   & 1826  & 0.12377 & 8.63E-07 & F     & F     & F     & F     & 39    & 49    & 0.037398 & 7.66E-07 \bigstrut\\
			\cline{2-14}          & $x_0^5$ & 857   & 2987  & 0.50011 & 8.71E-07 & F     & F     & F     & F     & 42    & 55    & 0.047374 & 9.82E-07 \bigstrut\\
			\cline{2-14}          & $x_0^6$ & 759   & 2656  & 0.53768 & 8.85E-07 & F     & F     & F     & F     & 32    & 48    & 0.037981 & 8.13E-07 \bigstrut\\
\hline
		\end{tabular}%
	}
	\label{tab6}%
\end{table}%

\begin{table}[htbp]
	\centering
	\caption{Test results for the \texttt{hybrid+}, \texttt{dfrmil} and \texttt{dfprp} methods on problem 7}
	\resizebox{\columnwidth}{!}{
		\begin{tabular}{|c|r|c|c|c|c|c|c|c|c|c|c|c|c|}
			\hline
			\multirow{2}[4]{*}{Dimension} & \multicolumn{1}{c|}{\multirow{2}[4]{*}{$x_0$}} & \multicolumn{4}{c|}{dfrmil}   & \multicolumn{4}{c|}{dfprp}    & \multicolumn{4}{c|}{hybrid+} \bigstrut\\
			\cline{3-14}          & \multicolumn{1}{c|}{} & \#iter & \#fval & time  & Fnorm & \#iter & \#fval & time  & Fnorm & \#iter & \#fval & time  & Fnorm \bigstrut\\
			\hline
			\multirow{10}[20]{*}{1000} & $x_0^1$ & F     & F     & F     & F     & 25    & 25    & 0.021629 & 4.58E-07 & 28    & 54    & 0.043669 & 1.28E-08 \bigstrut\\
			\cline{2-14}          & $x_0^2$ & F     & F     & F     & F     & 21    & 21    & 0.014976 & 3.06E-07 & 24    & 45    & 0.024179 & 1.48E-07 \bigstrut\\
			\cline{2-14}          & $x_0^3$ & 232   & 379   & 0.13626 & 9.21E-07 & 25    & 25    & 0.012154 & 7.59E-07 & 8     & 8     & 0.006353 & 5.55E-07 \bigstrut\\
			\cline{2-14}          & $x_0^4$ & 281   & 536   & 0.16753 & 4.36E-07 & 31    & 31    & 0.01693 & 6.01E-07 & 27    & 28    & 0.024567 & 8.32E-07 \bigstrut\\
			\cline{2-14}          & $x_0^5$ & 380   & 629   & 0.095196 & 9.49E-07 & 30    & 30    & 0.01731 & 9.66E-07 & 32    & 40    & 0.045519 & 5.07E-07 \bigstrut\\
			\cline{2-14}          & $x_0^6$ & 195   & 287   & 0.10035 & 8.76E-07 & 27    & 27    & 0.022926 & 6.09E-07 & 18    & 18    & 0.025864 & 3.47E-07 \bigstrut\\
\hline
		\end{tabular}%
	}
	\label{tab7}%
\end{table}%

\begin{table}[htbp]
	\centering
	\caption{Test results for the \texttt{hybrid+}, \texttt{dfrmil} and \texttt{dfprp} methods on problem 8}
	\resizebox{\columnwidth}{!}{
		\begin{tabular}{|c|r|c|c|c|c|c|c|c|c|c|c|c|c|}
			\hline
			\multirow{2}[4]{*}{Dimension} & \multicolumn{1}{c|}{\multirow{2}[4]{*}{$x_0$}} & \multicolumn{4}{c|}{dfrmil}   & \multicolumn{4}{c|}{dfprp}    & \multicolumn{4}{c|}{hybrid+} \bigstrut\\
			\cline{3-14}          & \multicolumn{1}{c|}{} & \#iter & \#fval & time  & Fnorm & \#iter & \#fval & time  & Fnorm & \#iter & \#fval & time  & Fnorm \bigstrut\\
			\hline
			\multirow{9}[18]{*}{1000} & $x_0^1$ & F     & F     & F     & F     & 33    & 33    & 0.062222 & 6.04E-07 & 3     & 3     & 0.01581 & 3.19E-09 \bigstrut\\
			\cline{2-14}          & $x_0^2$ & 2     & 2     & 0.005009 & 3.16E-08 & 28    & 28    & 0.035583 & 8.87E-07 & 2     & 2     & 0.00654 & 3.16E-08 \bigstrut\\
			\cline{2-14}          & $x_0^3$ & 3     & 3     & 0.003307 & 3.05E-07 & 25    & 25    & 0.026045 & 8.05E-07 & 3     & 3     & 0.006138 & 1.56E-07 \bigstrut\\
			\cline{2-14}          & $x_0^4$ & 4     & 4     & 0.007079 & 2.31E-08 & 32    & 32    & 0.047148 & 5.98E-07 & 4     & 4     & 0.006411 & 3.27E-09 \bigstrut\\
			\cline{2-14}          & $x_0^5$ & 4     & 4     & 0.021576 & 2.34E-08 & 32    & 32    & 0.039794 & 5.98E-07 & 4     & 4     & 0.008313 & 9.43E-12 \bigstrut\\
			\cline{2-14}          & $x_0^6$ & 4     & 4     & 0.004315 & 5.26E-13 & 27    & 27    & 0.045387 & 6.13E-07 & 4     & 4     & 0.010258 & 1.10E-11 \bigstrut\\
\hline
		\end{tabular}%
	}
	\label{tab8}%
\end{table}%

\begin{table}[htbp]
	\centering
	\caption{Test results for the \texttt{hybrid+}, \texttt{dfrmil} and \texttt{dfprp} methods on problem 9}
	\resizebox{\columnwidth}{!}{
		\begin{tabular}{|c|r|c|c|c|c|c|c|c|c|c|c|c|c|}
			\hline
			\multirow{2}[4]{*}{Dimension} & \multicolumn{1}{c|}{\multirow{2}[4]{*}{$x_0$}} & \multicolumn{4}{c|}{dfrmil}   & \multicolumn{4}{c|}{dfprp}    & \multicolumn{4}{c|}{hybrid+} \bigstrut\\
			\cline{3-14}          & \multicolumn{1}{c|}{} & \#iter & \#fval & time  & Fnorm & \#iter & \#fval & time  & Fnorm & \#iter & \#fval & time  & Fnorm \bigstrut\\
			\hline
			\multirow{10}[20]{*}{1000} & $x_0^1$ & 605   & 1027  & 0.38211 & 9.29E-07 & 15    & 15    & 0.047734 & 5.42E-07 & 6     & 7     & 0.21504 & 9.60E-08 \bigstrut\\
			\cline{2-14}          & $x_0^2$ & 583   & 990   & 0.14672 & 9.83E-07 & 15    & 15    & 0.05347 & 5.46E-07 & 6     & 7     & 0.056606 & 1.68E-11 \bigstrut\\
			\cline{2-14}          & $x_0^3$ & F     & F     & F     & F     & 15    & 15    & 0.011386 & 6.98E-07 & 6     & 7     & 0.11148 & 5.87E-10 \bigstrut\\
			\cline{2-14}          & $x_0^4$ & F     & F     & F     & F     & 15    & 15    & 0.011298 & 3.44E-07 & 6     & 7     & 0.008645 & 1.05E-08 \bigstrut\\
			\cline{2-14}          & $x_0^5$ & F     & F     & F     & F     & 15    & 15    & 0.019556 & 3.44E-07 & 6     & 7     & 0.022046 & 1.05E-08 \bigstrut\\
			\cline{2-14}          & $x_0^6$ & F     & F     & F     & F     & 15    & 15    & 0.014242 & 6.90E-07 & 6     & 7     & 0.020666 & 4.13E-07 \bigstrut\\
\hline
		\end{tabular}%
	}
	\label{tab9}%
\end{table}%

\begin{table}[htbp]
	\centering
	\caption{Test results for the \texttt{hybrid+}, \texttt{dfrmil} and \texttt{dfprp} methods on problem 10}
	\resizebox{\columnwidth}{!}{
		\begin{tabular}{|c|r|c|c|c|c|c|c|c|c|c|c|c|c|}
			\hline
			\multirow{2}[4]{*}{Dimension} & \multicolumn{1}{c|}{\multirow{2}[4]{*}{$x_0$}} & \multicolumn{4}{c|}{dfrmil}   & \multicolumn{4}{c|}{dfprp}    & \multicolumn{4}{c|}{hybrid+} \bigstrut\\
			\cline{3-14}          & \multicolumn{1}{c|}{} & \#iter & \#fval & time  & Fnorm & \#iter & \#fval & time  & Fnorm & \#iter & \#fval & time  & Fnorm \bigstrut\\
			\hline
			\multirow{10}[20]{*}{1000} & $x_0^1$ & F     & F     & F     & F     & 33    & 33    & 0.028738 & 6.65E-07 & 6     & 7     & 0.01399 & 1.95E-12 \bigstrut\\
			\cline{2-14}          & $x_0^2$ & F     & F     & F     & F     & 28    & 28    & 0.011566 & 8.89E-07 & 3     & 3     & 0.014673 & 1.61E-10 \bigstrut\\
			\cline{2-14}          & $x_0^3$ & 27    & 42    & 0.019897 & 3.44E-07 & 25    & 25    & 0.013694 & 8.33E-07 & 6     & 6     & 0.015631 & 8.13E-13 \bigstrut\\
			\cline{2-14}          & $x_0^4$ & 113   & 232   & 0.030826 & 4.06E-07 & 32    & 32    & 0.007726 & 6.49E-07 & 8     & 9     & 0.016942 & 2.47E-08 \bigstrut\\
			\cline{2-14}          & $x_0^5$ & 161   & 293   & 0.082147 & 4.52E-07 & 32    & 32    & 0.017012 & 6.49E-07 & 8     & 9     & 0.022529 & 2.47E-08 \bigstrut\\
			\cline{2-14}          & $x_0^6$ & 19    & 41    & 0.018665 & 5.17E-07 & 27    & 27    & 0.020302 & 6.67E-07 & 6     & 6     & 0.010383 & 2.59E-08 \bigstrut\\
\hline
		\end{tabular}%
		\label{tab10}%
	}
\end{table}%
\end{document}